\documentclass[final,5p,times,twocolumn]{elsarticle}

\newcommand\blfootnote[1]{%
  \begingroup
  \renewcommand\thefootnote{}\footnote{#1}%
  \addtocounter{footnote}{-1}%
  \endgroup
}

\usepackage{ragged2e}
\usepackage{etoolbox}
\usepackage{lipsum}

\usepackage[english]{babel}
\usepackage{tikz}
\usepackage[T1]{fontenc}
\usepackage[utf8]{inputenc}
\usepackage{ae}
\usepackage{amsmath}
\usepackage{amssymb} 
\usepackage{latexsym}
\usepackage{graphicx,epsfig,psfrag}
\usepackage{amsmath}
\usepackage{rotate}
\usepackage{amsfonts}
\usepackage{bbold}

\newtheorem{proposition}{Proposition}
\newtheorem{example}{Example}

\newtheorem{theorem}{Theorem}
\newtheorem{lemma}{Lemma}

\newtheorem{corollary}{Corollary}
\newcommand {\be}{\begin{equation}}
\newcommand {\ee}{\end{equation}}

\usepackage{subfigure}

\usepackage{tikz}
\usepgflibrary {shadings} 
\usetikzlibrary{3d,shadings}

\usepackage{accents}
\newcommand{\ubar}[1]{\underaccent{\bar}{#1}}

\def\ls#1{{\color{black}#1}}

\journal{Operations Research Letters}

\begin{document}

\begin{frontmatter}



\title{Cascading Failures in the Global Financial System: A Dynamical  Model}


\author[1]{Leonardo~Stella}
\author[2,3]{Dario~Bauso}
\author[4]{Franco~Blanchini}
\author[5]{Patrizio~Colaneri}

\affiliation[1]{organization={School of Computer Science, University of Birmingham},
            country={United Kingdom}}

\affiliation[2]{organization={Jan C. Willems Center for Systems and Control, ENTEG, University of Groningen},
            country={The Netherlands}}
            
\affiliation[3]{organization={Dipartimento di Ingegneria, University of Palermo},
            country={Italy}}
            
\affiliation[4]{organization={Dipartimento di Matematica e Informatica, Universit\`a degli Studi di Udine},
            country={Italy}}
            
\affiliation[5]{organization={Dipartimento di Elettronica, Informazione e Bioingeneria, and IEIIT-CNR, Politecnico di Milano},
            country={Italy}}

            

\begin{abstract}
\ls{In this paper, we propose a dynamical model to capture cascading failures among interconnected organizations in the global financial system and develop a framework to investigate under which conditions organizations remain healthy. The contribution of this paper is threefold: i) we develop a dynamical model that describes the time evolution of the organizations' equity values given nonequilibrium initial conditions; ii) we characterize the equilibria for this model; and iii) we provide a computational method to anticipate potential propagation of failures.}
\end{abstract}




\end{frontmatter}


\section{Introduction}\label{sec:intro} \blfootnote{\ls{Corresponding author: L. Stella; address: University of Birmingham, Edgbaston, Birmingham, B15 2TT; e-mail: l.stella@bham.ac.uk.}} 
In the wake of recent events concerning the collapses of Silicon Valley Bank and Credit Suisse (CS), the focus of this paper is to investigate the propagation of failures in financial systems. The  current global financial system is the resultant of a large number of financial interdependencies among governments, banks, firms, smaller and larger companies, private citizens, etc. In the same spirit of the related literature, we make use of the term organization in a broad sense including all these entities and individuals. These organizations hold each other's shares, debts and obligations in variable proportions. As a result, when a failure occurs, this can propagate through the network of interdependencies bringing other organizations to bankruptcy. \ls{Failures can take the form of bankruptcies, defaults, and other insolvencies.} Indeed, cascading defaults and failures account for one of the highest risks for the global financial system, let aside those institutions that are considered too big to fail, e.g., central banks. A slightly less recent example, but equally prominent, can be found in the interventions put together by the European Commission to save Greece and Spain from default following the historic quote ``whatever it takes'' by ECB President Mario Draghi (July 23rd, 2012) \cite{Elliott_2014}.

In this paper, we study the role of cascading failures among organizations linked through a network of financial interdependencies in the global financial system. Our aim is to develop a model that describes the risks associated with the propagation of failures in the network as well as the design of effective responses to mitigate the impact of financial contagion. Indeed, in the proposed model we highlight three relevant aspects: i) the interdependencies in a financial system through cross-holdings of shares or other liabilities; ii) the market price of assets owned by each organization; iii) and a failure cost incurred by each organization. Indeed, when the value of a financial organization attains a value that is below a failure threshold, additional losses propagate through the network leading to a cascade of failures.

\textit{Related Works}. The first structural framework to study the propagation of shocks in inter-bank lending was originally proposed in a pioneering work by Eisenberg and Noe in 2001 \cite{Eisenberg_2001}. The main contribution of that work is the introduction of a model that captures the contagion from individual organizations to other organizations in an inter-bank lending network. The contagion occurs at individual nodes and propagates in the network, leading to new equilibrium points representing the agreed mutual payments.

Their model illustrated how shocks to individual organizations can propagate through inter-bank lending networks. Indeed, contagion develops instantaneously, bringing the network to a new equilibrium on an agreed set of mutual payments.  Subsequently, there has been a substantial body of work analyzing and generalizing this framework. For example, the authors in \cite{Carvalho_2008} and \cite{Acemoglu_2012} studied the way in which the structure of network graphs, such as hubs, sparsity, and asymmetry structure, influences the shock propagation and the magnitude of the aggregate fluctuation. Their study provides insights on the optimal structure for inter-bank lending networks. Their model can accommodate a variety of settings. For production networks, the model represents the input-output relationship and determines the output equilibrium \cite{Acemoglu_2012}, whereas for financial systems, it calculates the clearing loan repayments, involving the systemic risk of default cascade~\cite{Eisenberg_2001}.

Later, the preliminary research proposed by Eisenberg and Noe was extended in several directions. A body of literature dating back to the work by Elsinger \cite{Elsinger_2009} and then followed by Elliott \emph{et al.} \cite{Elliott_2014}, Rogers and Veraart \cite{Rogers_2013}, and Glasserman and Young \cite{Glasserman_2015} considered bankruptcy costs and their impact onto the financial system. As a consequence of these costs, financial organizations can in turn fail and drag other organizations to bankruptcy. Simultaneously, cross-holdings were considered by Elsinger \cite{Elsinger_2009}, Elliott \emph{et al.} \cite{Elliott_2014}, Fischer \cite{Fischer_2014} and Karl and Fischer \cite{Karl_2014}. An important aspect in many of these works is that cross-holdings inflate the value of the financial system and thus the net value of each organization needs to be adjusted by a factor that preserves the real value in the system \cite{Brioschi_1989}. The work by Weber and Weske considers both these aspects and integrates them into a system that is able to capture fire sales as well \cite{Weber_2017}.

In particular, the work by Elliott \emph{et al.} highlighted the fact that in the current highly interconnected financial system, where banks and other institutions are linked via a network of mutual liabilities, a financial shock in one or few nodes of the network may hinder the possibility for these nodes to fulfill their obligations towards other nodes, and therefore provoke default \cite{Elliott_2014}. A recent work by Birge \cite{Birge_2021} investigates an inverse optimisation approach based on the decisions from national debt cross-holdings to address the propagation and extent of failures in the network.

However, the common assumption that all payments are simultaneous is quite unrealistic. For this reason, several recent works, e.g.,  see \cite{Acemoglu_2015, Cabrales_2017, Chen_2021, Glasserman_2016}, propose time-dynamic extensions of this model. The work by Calafiore \emph{et al.} considers the problem of reducing the financial contagion by introducing some targeted interventions that can mitigate the cascaded failure effects. They consider a multi-step dynamic model of clearing payments and introduce an external control term that represents corrective cash injections made by a ruling authority \cite{Calafiore_2022}. Similarly, a case study on the Korean financial system is proposed by Ahn and Kim where the authors study the interventions in the form of liquidity injection into the financial system under economic shocks \cite{Ahn_2019}. Finally, a recent work by Ramirez \emph{et al.} investigated a stochastic discrete-time model where the mean and covariance error are studied with focus on the steady-state solution \cite{Ramirez_2022}.

\textit{Contribution}. The contribution of this work is threefold. Firstly, we introduce the formulation of a dynamical model for cascading failures in financial systems. \ls{This model is novel with respect to the literature, and in particular to the work by Elliott \emph{et al.}, as in the following:
\begin{itemize}
\item Our model captures the transient response, allowing us to study the market response to disturbances, and uncertainty in the form of initial conditions not already at an equilibrium.
\item Our model can predict the future evolution of the market, allowing us to characterize the equilibria and study local stability.
\item Finally, it allows for the study of sensitivity with respect to the parameters. Moreover, in the case of time-varying parameters, e.g., the prices of assets, our model is able to accommodate for fluctuations and convergence to a stable trajectory.
\end{itemize}
The second} contribution of this paper is the stability analysis of the equilibrium points of the proposed system. In particular, we show the existence of these equilibria, their uniqueness and provide an explicit expression for them. \ls{The third contribution is} a computational method via sign-space iteration that allows us to compute the attractive equilibrium point for given initial conditions.

The paper is organized as follows. First, we introduce the notation. In Section~\ref{sec:problem}, we develop the networked model. In Section~\ref{sec:stability}, we investigate the existence, uniqueness and stability of the equilibrium points of our system. In Section~\ref{sec:sign}, we illustrate the computational algorithm. Finally, in Section~\ref{sec:conc}, we discuss concluding remarks and future directions.

\medskip
\noindent \textbf{Notation}. 
The symbols $\mathbb{0}_n$ and $\mathbb{1}_n$ denote the $n$-dimensional column vector with all entries equal to 0 and to 1, respectively. The identity matrix of order $n$ is denoted by $I_n$. Let $J^{[k]} := {\rm diag} (1 - 2\phi^{[k]})$, where vector $\phi^{[k]}$ represents the integer $k$ in binary representation; we denote the generic orthant $k$ by $\mathcal X^k$, namely, $\mathcal X^k := \{x \in \mathbb R^n | J^{[k]}x \ge 0 \}$. Given a generic vector $V \in \mathbb R^n$, let the operator $y = \phi(V)$ be such that the $i$th component \ls{$y_i = 0$ if $V_i \ge 0$ or $y_i = 1$}, otherwise. The notation $V \ge 0$ for a generic vector $V$ or $M \ge 0$ for a generic matrix $M$ is to be intended elementwise.

A square real matrix $M \in \mathbb R^{n\times n}$ is said to be \emph{Metzler} if its off-diagonal entries are nonnegative, namely, $M_{i,j} \ge 0$, $i \neq j$. Every Metzler matrix $M$ has a real dominant eigenvalue $\lambda_F(M)$, which is referred to as \emph{Frobenius eigenvalue}. The corresponding left and right vectors associated with $\lambda_F(M)$ are referred to as left and right \emph{Frobenius eigenvectors}, respectively. A square real matrix $M$ is said to be \emph{Hurwitz} if all its eigenvalues lie in the open left half plane. A square matrix is said to be \ls{\emph{Schur}} if all its entries are real and its eigenvalues have absolute value less than one \cite{Farina_2000}.

\section{Problem Formulation}\label{sec:problem}
In this section, we introduce the model of a networked financial system, where a number of organizations are linked through financial interdependencies. To this aim, we consider a set of organizations $N = \{1, \dots, n\}$. Each organization $i \in N$ is described by an equity value $V_i \in \mathbb R$, which represents the total values of its shares.
Organizations can invest in primitive assets, namely, mechanisms that generate income in the form of a net flow of cash over time. We consider a set of primitive assets $M = \{1, \dots, m\}$. We denote the market price of asset $k$ by $p_k$ and the share of the value of asset $k$ held by organization $i$ by $D_{ik} \ge 0$. Each organization can also hold shares of other organizations; for any pair of organizations $i,j \in N$, let $C_{ij} \ge 0$ be the fraction of organization $j$ owned by organization $i$.

The equity values of organizations can be determined by the following discrete-time dynamical model:
\begin{equation}\label{eq:model}
V(t+1)=CV(t)+Dp-B\phi(V(t)-\ubar V),
\end{equation}
where $t \in {\mathbb Z}^+$, $C$ is a nonnegative and nonsingular matrix where $C_{ii}=0$ and $\mathbb{1}_n^\top  C< \mathbb{1}_n^\top$ which means that the equity value of each organization held by other organizations cannot exceed the equity value of the organization itself,  $D$ is a positive matrix, $p$ a nonnull nonnegative vector, $B={\rm diag}(\beta)$ \ls{is a nonnegative diagonal matrix} with entries $\beta_i>0$, $i=1,2,\cdots, N$, $\ubar V$ is the vector of threshold values $\ubar V_i$ for all $i$ below which organization $i$ incurs a failure cost $\beta_i$ and $\phi(V-\ubar V)$ the vector of indicator functions taking value $1$ if $V_i<\ubar V_i$ and $0$ if $V_i\ge \ubar V_i$. The first term in~\eqref{eq:model} takes into account the cross-holdings, the second term describes the primitive assets held by each organization and the last term accounts for the discontinuous drop imposed by the cost of failure.

\ls{The main difference with the papers in the literature is that our model, namely, system~\eqref{eq:model}, captures the impact of the transient onto the steady-state market values. In order to emphasize this, we present the following example.} 

\begin{example}\label{ex1}
\ls{Consider system~(\ref{eq:model}) with $N = 2$ organisations, $M=2$ assets. The initial condition $V(0)$ is set to be random in $[2, 5]$. Let $C = [0 \; 0.025; \; 0.005 \; 0]$, $D = [0.05 \; 0.05; \; 0.05 \; 0.05]$, $\beta = \mathbb 1_{2}$, and $\ubar V = 1.5 \; \mathbb 1_{20}$.}

\ls{Now, let us consider two scenarios and let the asset price be a time-varying signal $p(t)$ that changes over time. In the first scenario, $p(t)$ is set to 20 at the start of the simulation for both assets, $M=1,2$; a perturbation of one time instant in length occurs at $t=4$, making the price drop to 14.9. Likewise, in the second scenario, $p(t)$ is set to 20 at the start of the simulation for both assets, $M=1,2$; the size of the perturbation is the same as before, namely changing the value of $p(t)$ to 14.9, but the time window in which the perturbation occurs spans many time instants, namely, from $t=4$ to $t=20$. Figure~\ref{fig:pt} depicts the two scenarios (top-left and bottom-left) and the corresponding $p(t)$ (top-right and bottom-right). Figure~\ref{fig:pt} (top) shows the time evolution of the system, where the dashed red line represents the threshold: in this scenario, both companies remain healthy. Figure~\ref{fig:pt} (bottom) shows the situation where a longer-lasting perturbation affects the dynamics leading to an equilibrium point where one company fails.}

\ls{\textit{Remark}. Despite the simplicity of this toy problem, the example shows the ability of our model to capture the impact of the transient response to the system dynamics and, thus, the convergence to other potential equilibria, even though the final value of $p(t)$ is the same in both scenarios. This nonlinear behavior marks the difference with what was previously investigated in the literature and allows us to study perturbations of the market rather than just the system state at an equilibrium.}
\begin{figure}[t]
    \centering
    \includegraphics[width=0.5\textwidth]{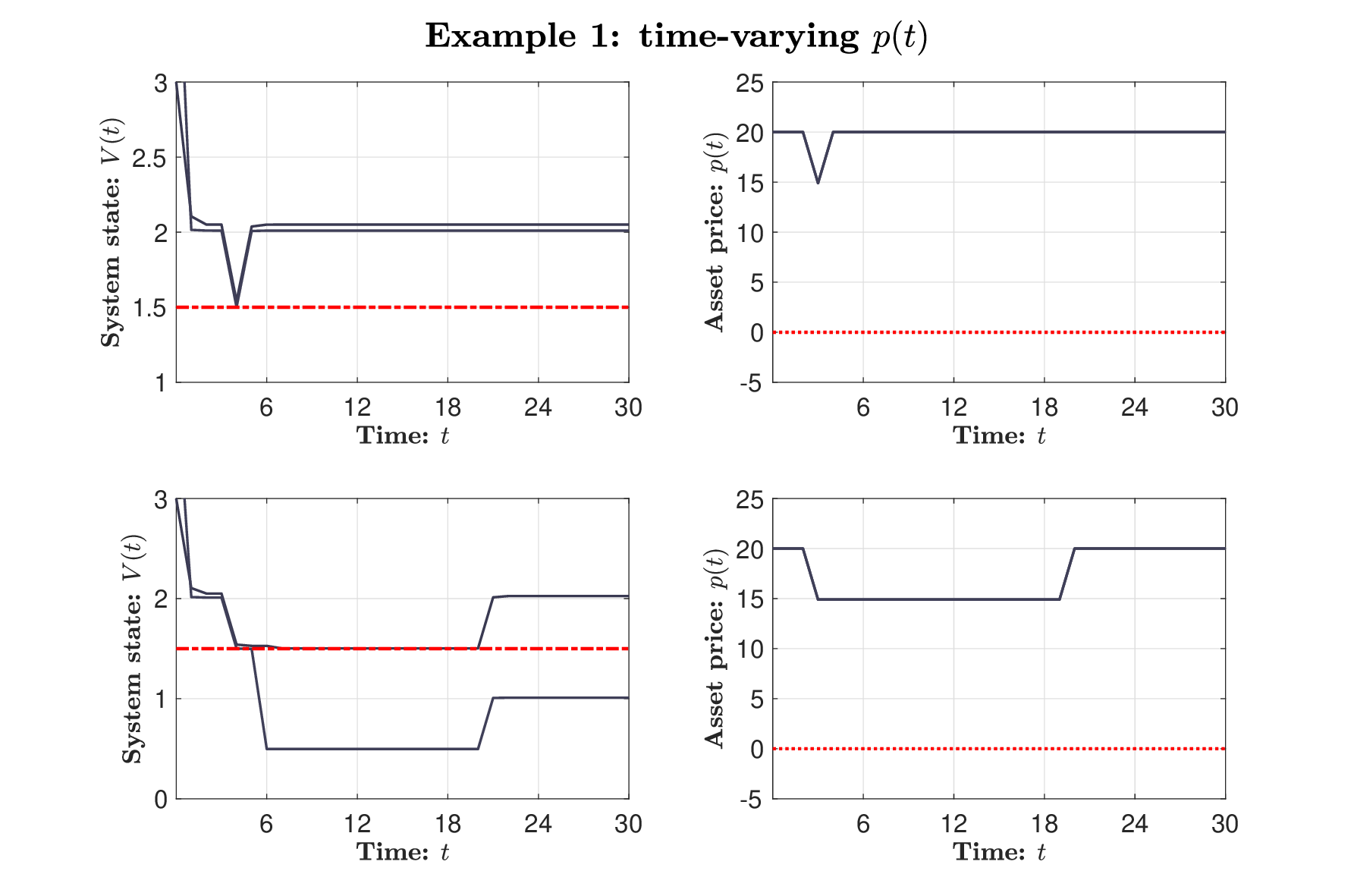}
    \caption{\ls{Example 1: given identical initial conditions and parameters, we consider two scenarios where $p(t)$ is a time-varying signal with same size of perturbation but different duration. In the first scenario, both companies remain healthy (top), whereas, due to the transient, one company fails in the second scenario (bottom).}}
    \label{fig:pt}
\end{figure}
\end{example}

\section{Characterization of the Equilibria}\label{sec:stability}
In this section, we study the equilibria of system~(\ref{eq:model}). From the condition that $0\le B\phi(V(t)-\ubar V)\le B$, and recalling that $C$ is nonnegative we derive the following preliminary result. 

\medskip\noindent
\begin{theorem}\label{th.1}
$V(t)\ge 0,  \forall t\ge 0$ and $V(0) \ge \mathbb 0_n$ if and only if 
\begin{equation}\label{eq.2}
Dp-\beta\ge 0.
\end{equation}  \hfill $\square$
\end{theorem}
Under condition (\ref{eq.2}), system (\ref{eq:model}) is a positive nonlinear switched system since vector $\phi(V(t)-\ubar V)$ can take a finite number of values $\phi^{[k]}$, with $k=0,1,2,\cdots, 2^n-1$. For instance, with $n=3$ we have:
\begin{equation}\nonumber
\begin{array}{c}
\phi^{[0]}= \mathbb{0}_n, \phi^{[1]}_3=\left[\begin{array}{c} 0 \\ 0 \\ 1\end{array}\right], \phi^{[2]}=\left[\begin{array}{c}  0 \\ 1 \\ 0\end{array}\right], \phi^{[3]}=\left[\begin{array}{c}  0 \\ 1 \\ 1\end{array}\right], \\
\phi^{[4]}=\left[\begin{array}{c} 1 \\ 0 \\ 0\end{array}\right], \phi^{[5]}=\left[\begin{array}{c} 1 \\ 0 \\ 1\end{array}\right], \phi^{[6]}=\left[\begin{array}{c} 1 \\ 1 \\ 0\end{array}\right], \phi^{[7]}= \mathbb{1}_n.
\end{array}
\end{equation}

As such, system~(\ref{eq:model}) may possess at most $2^n$ equilibria in total. The equilibria in orthant $k$, denoted by $\overline V^{[k]}$ and characterized by the index $k$ is given by 
\begin{equation}
\label{eq.3}
\overline V^{[k]}=(I_n-C)^{-1}(Dp-B\phi^{[k]}), \quad s.t. \quad \phi(\overline V^{[k]}-\ubar V)=\phi^{[k]}.
\end{equation} 
Note that $V=0$ cannot be an equilibrium of the system since $Dp > 0$ and that, if (\ref{eq.2}) holds, $\overline V^{[k]} > 0$. In the $k$th orthant the difference $Y^{[k]}(t)=V(t) - \overline V^{[k]}$ follows the autonomous dynamics 
\begin{equation}\label{eq.4}
Y^{[k]}(t+1)=CY^{[k]}(t).
\end{equation}
Since $C$ is nonnegative with $\mathbb{1}_n^\top  C< \mathbb{1}_n^\top$, it turns out that $C$ is Schur-stable. Therefore, the following theorem can be stated.


\medskip\noindent
\begin{theorem}\label{th.2}
Any equilibrium $\overline V^{[k]}$ which is in the interior of the $k$th orthant $\mathcal X^k$ is locally asymptotically stable. \hfill $\square$
\end{theorem}

\textit{Remark}. Note that there could be equilibria on the discontinuity points, 
but these are fragile (unstable) and are not considered.

\begin{figure}[b]
    \centering
    \includegraphics[width=0.49\textwidth]{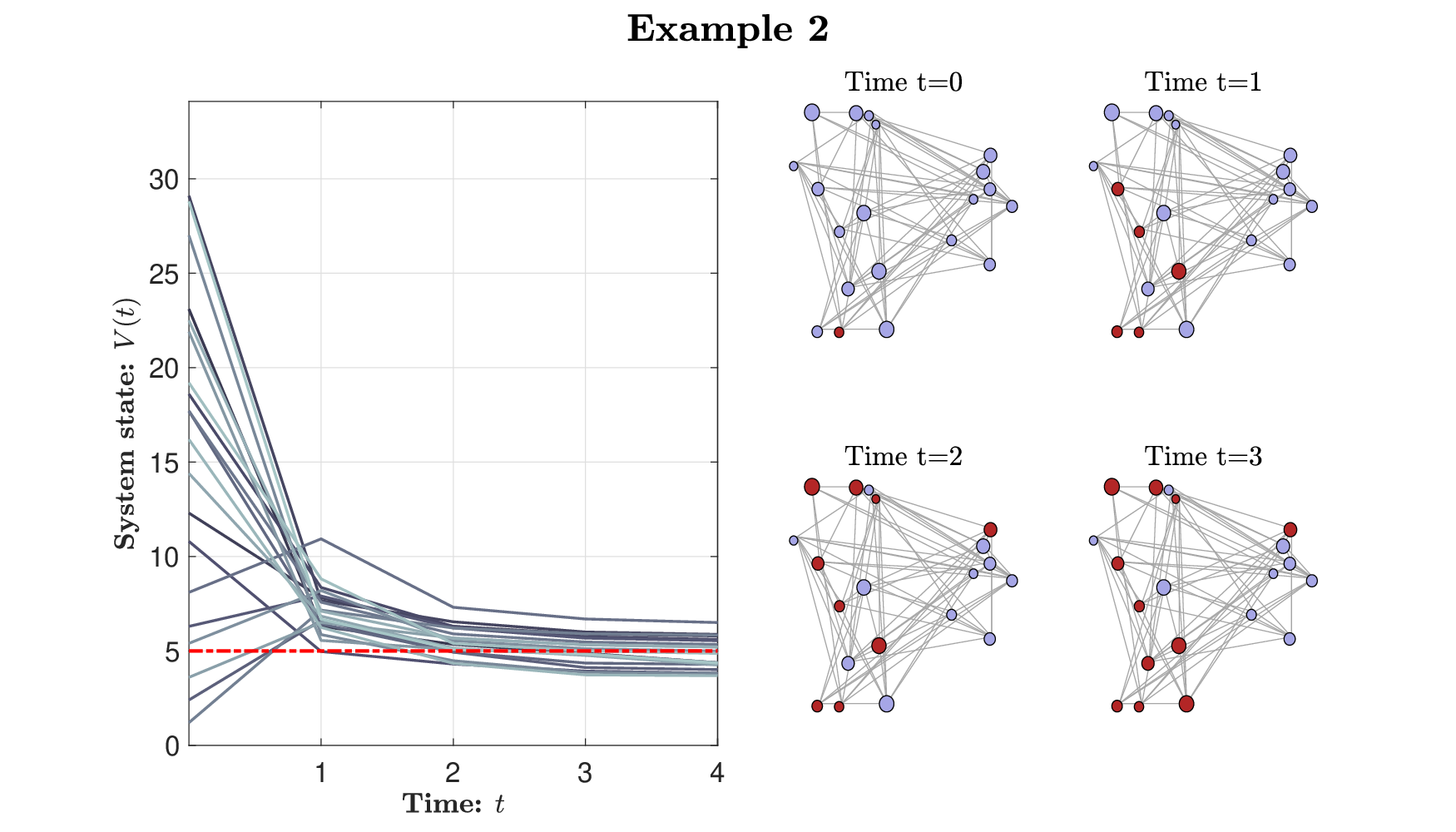}
    \caption{Example 1: since condition \eqref{eq.2} is satisfied, $V(t) \ge 0, \forall t \ge 0$ (left); network topology in the first four time instants (right).}
    \label{fig:ex2}
\end{figure}
\begin{example}\label{ex2}
Consider system~(\ref{eq:model}) with $N = 20$ organisations, $M=10$ assets. The initial condition $V(0)$ is set to be random in $[0, 30]$. Let $C$ be set to random values in $[0, 0.01]$ such that $C_{ii}=0$ and $\mathbb{1}_n^\top  C< \mathbb{1}_n^\top$. Finally, let
\begin{eqnarray}\nonumber
& D = 0.05 \; \mathbb 1_{20} \mathbb 1^\top_{10} \quad p = 10 \; \mathbb 1^\top_{10}, \\ 
\nonumber & \beta = \mathbb 1_{20}, \quad \ubar V = 10 \; \mathbb 1_{20}.
\end{eqnarray}
It is straightforward to see that $Dp-\beta = 4 \mathbb 1^\top_{20} \ge 0$. Therefore, in accordance to Theorem \ref{th.1}, the values of all companies remain positive, namely, $V(t) \ge 0, \forall t \ge 0$. Figure~\ref{fig:ex2} depicts this scenario. Figure~\ref{fig:ex2} (left) shows the time evolution of the system, where the dashed red line represents the threshold. Figure~\ref{fig:ex2} (right) shows the network topology in first four instants, where companies are indicated by colored nodes and edges indicate the cross-holdings between companies: the companies whose values are above the threshold are indicated in \ls{blue}, and in \ls{red, otherwise}.

\ls{\textit{Remark}. This example shows that if the condition in Theorem~\ref{th.1} holds true, the equity values of all organizations will remain positive at all time. This provides an important addition to previous studies on this topic, as we are able to predict the behavior of the system at every time instant.}
\end{example}

\medskip
We now turn our attention to the existence and uniqueness of the equilibrium points in orthants $0$ and $2^n-1$, which we henceforth refer to as \emph{positive} and \emph{negative} equilibrium points, respectively. To this aim, consider:
 \begin{equation}\nonumber
\begin{array}{c}
V(t+1) = CV(t) + Dp - B \phi(V(t) - \ubar V), \\
x(t) = V(t) - \ubar V.
\end{array}
\end{equation}
The above system can be rewritten as
\begin{equation}\label{eq:sysxr}
\begin{array}{c}
x(t+1) := Cx(t) + r - B \phi(x(t)), \\
r := (C - I_n) \ubar V + Dp.
\end{array}
\end{equation}
The above is a monotone system since $\phi(y) \ge \phi(x)$ if $y \le x$. We can now prove the following theorem.

\medskip
\begin{theorem}\label{th.negpos}
Consider system~\eqref{eq:sysxr}. In each open orthant $\mathcal X^k$, there exists at most one equilibrium. Furthermore, the following points hold true:
\begin{enumerate}
\item There exists an equilibrium point $\bar x \ge 0$ if and only if $(I_n - C)^{-1}r \ge 0$.
\item If  $(I_n - C)^{-1}(r - \beta) \ge 0$, then there exists an equilibrium point $\bar x \ge 0$ and it is the unique equilibrium.
\item There exists an equilibrium point $\bar x < 0$ if and only if $(I_n - C)^{-1}(r - \beta) < 0$.
\item If  $(I_n - C)^{-1}r < 0$, then there exists an equilibrium point $\bar x < 0$ and it is the unique equilibrium.
\end{enumerate} \hfill $\square$
\end{theorem}
\textit{Proof}. First, let us prove the first statement, namely, if an equilibrium exists in orthant $k$, it is unique. Let 
$$\bar x^{[k]} = (I_n - C)^{-1}(r - B \phi^{[k]}) \in \mathcal X^k$$
be the generic equilibrium point in the $k$th orthant. By contradiction, let us assume that a second equilibrium point exists in the same orthant. It is straightforward to see that the calculation with a given $\phi^{[k]}$ would produce the same equilibrium point. \ls{Note that in the rest of the proof, since $C$ is Schur, then $(C - I_n)$ is Hurwitz and Metzler and, therefore, the inverse of the negative, namely, $(I_n - C)^{-1} \ge 0$, elementwise \cite{Farina_2000}.}\\
Let us now prove the rest point by point.
\begin{enumerate}
\item Let $(I_n - C)^{-1}r \ge 0$, then $\bar x = (I_n - C)^{-1}r \ge 0 \in \mathcal X^0$. Vice versa, assume that there exists a generic equilibrium $\bar x \ge 0$, then $\bar x \in \mathcal X^0$. Therefore, $\phi(\bar x) = 0$ and $(I_n - C)^{-1}r \ge 0$.
\item Let $(I_n - C)^{-1}(r - \beta) \ge 0$, then $(I_n - C)^{-1}r \ge (I_n - C)^{-1} \beta \ge 0$. It follows from the first point that there exists an equilibrium $\bar x \ge 0$. Moreover, assume there exists an equilibrium $\bar x^{[k]}$ in orthant $\mathcal X^k$, i.e., \ls{$\bar x^{[k]} = (I_n -C)^{-1}(r - B \phi(x^{[k]})) \ge (I_n -C)^{-1}(r - \beta) \ge 0$}. Then, the unique equilibrium is in orthant $\mathcal X^0$. 
\item Let $(I_n - C)^{-1}(r - \beta) < 0$, then $\bar x = (I_n - C)^{-1}(r - \beta) < 0 \in \mathcal X^{2^n-1}$. Vice versa, assume that there exists a generic equilibrium $\bar x < 0$, then $\bar x  \in \mathcal X^{2^n-1}$. Therefore, $\bar x = (I_n - C)^{-1}(r - \beta) < 0$.
\item Let $(I_n - C)^{-1}r \le (I_n - C)^{-1}(r - \beta) < 0$, then from point 3, there exists an equilibrium $\bar x^{[k]} < 0$. Moreover, assume there exists an equilibrium $\bar x^{[k]}$ in orthant $\mathcal X^k$, i.e., $\bar x^{[k]} = (I_n -C)^{-1}r - B \phi(x^{[k]}) \le (I_n -C)^{-1}r < 0$. Then, the unique equilibrium is in orthant $\mathcal X^{2^n-1}$. 
\end{enumerate}
This concludes our proof. \hfill $\blacksquare$ \\

\begin{example}\label{ex3}
Consider system~(\ref{eq:sysxr}) with $N = 20$ organisations, $M=10$ assets. The initial condition $x(0)$ is set to be random in $[0, 30]$. Let $C$ be set to random values in $[0, 0.01]$ such that $C_{ii}=0$ and $\mathbb{1}_n^\top  C< \mathbb{1}_n^\top$. We provide two sets of simulations. Table~\ref{t:varying} includes all the other parameters for each simulation.

\begin{table} [h!]
\caption{Set of parameters for each simulation.}
\begin{center}
   \begin{tabular}{|c|c|c|c|c|}
  \hline
Simulation & $D$ & $p$ & $\beta$  & $\ubar V$ \\ \hline \hline
   I & $0.06 \; \; \mathbb 1_{20} \mathbb 1^\top_{10}$ & $10 \; \mathbb 1^\top_{20}$ & $\mathbb 1_{20}$ & $10 \; \mathbb 1_{20}$ \\ \hline 
   II & $0.03 \; \; \mathbb 1_{20} \mathbb 1^\top_{10}$ & $10 \; \mathbb 1^\top_{20}$ & $\mathbb 1_{20}$ & $10 \; \mathbb 1_{20}$ \\ \hline
    \end{tabular}\\
  \end{center}
  \label{t:varying}
\end{table}

In the first set of simulations, the positive equilibrium, namely, $\bar x \ge 0$ exists and is unique. This is in accordance with condition 1 and condition 2 of Theorem \ref{th.negpos}. This can be seen in Fig.~\ref{fig:ex3} (top-left). Similarly, in the second set of simulations,  since the third and last conditions of Theorem \ref{th.negpos} hold true, the negative equilibrium point, i.e., $\bar x < 0$, exists and is unique. Figure~\ref{fig:ex3} (bottom-left) shows the second set of simulations. Figure~\ref{fig:ex3} (right) shows the network topology in the first and third instant for each set of simulations. Colors have the usual meaning.
\begin{figure}[t]
    \centering
    \includegraphics[width=0.49\textwidth]{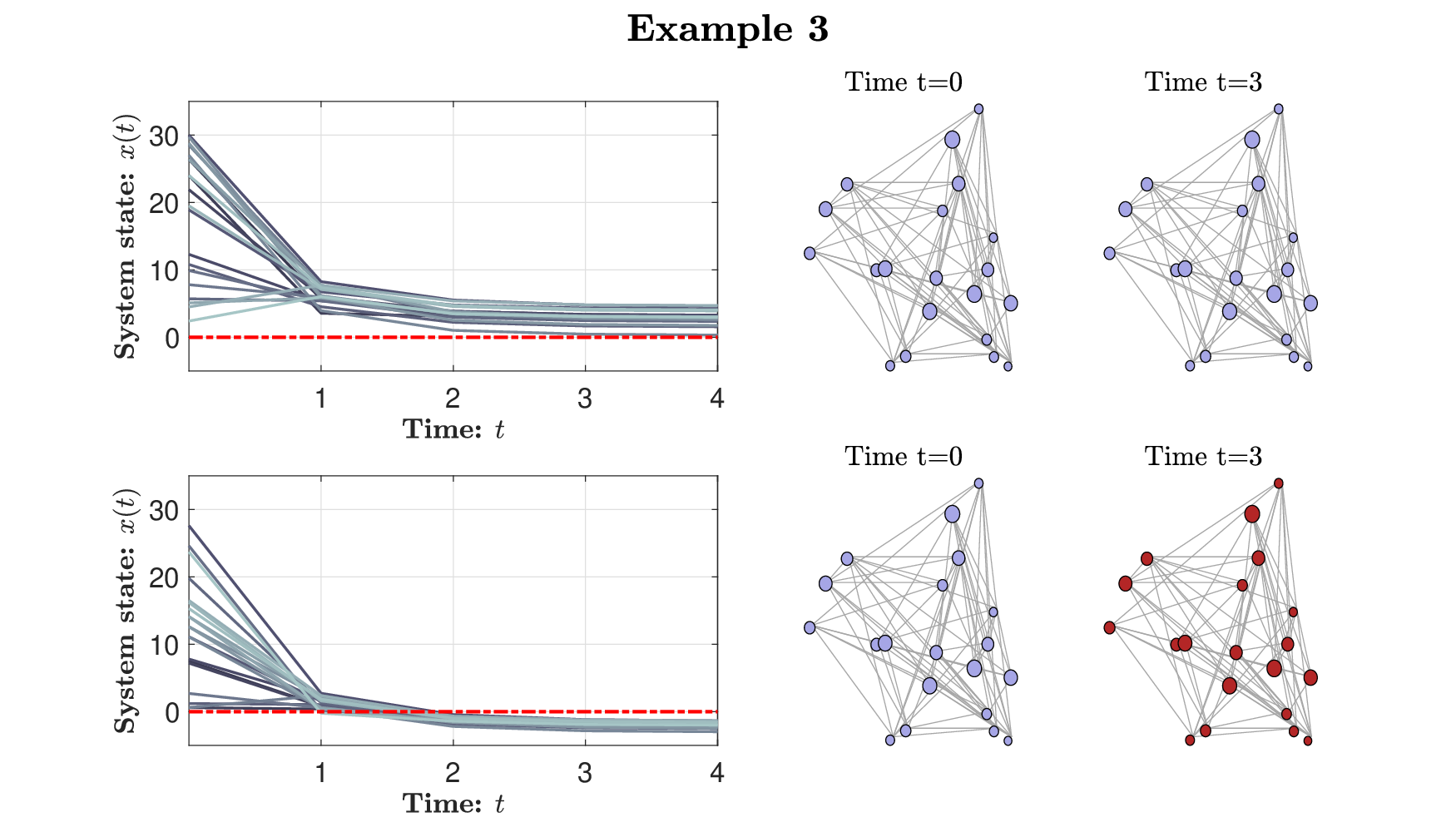}
    \caption{Example 2: the equilibrium point $\bar x \ge 0$ exists and is unique as condition 2 of Theorem \ref{th.negpos} holds true (top); similarly, since the 4-th condition of Theorem \ref{th.negpos} holds true, the equilibrium point $\bar x < 0$ exists and is unique (bottom).}
    \label{fig:ex3}
\end{figure}

\ls{\textit{Remark}. A physical interpretation of this example follows. The translated variable $x$ allows us to study the propagation of failures over time. If the conditions in Theorem~\ref{th.negpos} hold true, we can show whether the system converges to the equilibrium point in which all organizations are healthy or to the equilibrium point where all fail.}
\end{example}

Now, we provide a sufficient condition that guarantees that no equilibrium point in the negative orthant exists with respect to a subgraph of the cross-sharing $C$. 
\begin{proposition}
Given a square principal submatrix of $C$, denoted by $\tilde C$, if the following holds:
\begin{equation}\label{eq:subgraph}
\ubar V_i < \frac{(Dp - \beta)_i} {1 -\lambda_F (\tilde C)}, \quad \forall i,
\end{equation}
where $\tilde C$ is a principal sparse subgraph of $C$, then there does not exist the negative equilibrium point, i.e., at least one organization remains healthy.
\end{proposition}
\textit{Proof}. Assume that equation~\eqref{eq:subgraph} holds true. Since $\lambda_F(\tilde C) \le \lambda_F(C)$ \cite{Farina_2000}, then $$\ubar V_i < \frac{(Dp - \beta)_i}{1-\lambda_F (C)}, \quad \forall i.$$
Let $x > 0$ be the left Frobenius eigenvector of $C$, i.e., $x^\top C = \lambda_F(C)x^\top$. Then,
$$x^\top \ubar V < \frac{x^\top(Dp-\beta)}{1-\lambda_F(C)} = x^\top (I_n-C)^{-1}(Dp-\beta),$$
so that, being $Dp+(C-I_n)\ubar V = r$, we have:
\begin{equation*}
\begin{array}{lll}
x^\top \ubar V < x^\top (I_n-C)^{-1}(r - (C- I_n)\ubar V -\beta), \\
x^\top \ubar V < x^\top \ubar V + x^\top (I_n-C)^{-1}(r -\beta), \\
x^\top (I_n-C)^{-1}(r -\beta) > 0.
\end{array}
\end{equation*}
The above implies $(I_n-C)^{-1}(r-\beta) \nless 0$. From point 3 of Theorem \ref{th.negpos}, then no equilibrium $\bar x < 0$ exists. \hfill $\blacksquare$ \\


\textit{Remark}. Condition \eqref{eq:subgraph} provides a relation among three main elements of the original system: the thresholds, the underlying topology and the external assets. Therefore, since it is desirable that the system does not converge to the negative equilibrium point, by violating this condition on $\ubar V$ we ensure that at least one company is healthy. 

\section{Sign-space Iteration}\label{sec:sign}
In this section, we analyze the behavior of the trajectories of financial organizations that are below and above the threshold. To this end, let us rewrite system~\eqref{eq:sysxr} in a more compact way as:
\begin{equation}\label{eq:x}
x(t+1) = C x(t) + \Psi(x(t)),
\end{equation}
where $\Psi(x) := r - B\phi(x)$ and, in particular, with a slight abuse of notation, the following 
$$
\Psi(x) = \Psi(\mbox{sign}(x)),~~~\Psi_k \in \{\psi_k^-,\psi_k^+\}
$$
depends on the sign of $x$, $\psi_k^- = r_k - \beta_k$ and $\psi_k^+ = r_k$ can both take positive and negative values. Here, the $\mbox{sign}(x)$ function is defined as:
$$\mbox{sign}(x) := 1 -2\phi(x) = \left\{ \begin{array}{ll} +1, \quad \mbox{if} \; x\ge 0 \\
-1, \quad \mbox{if} \; x < 0.
\end{array} \right.$$ 
Let $P=(I-C)^{-1}$. Then, an explicit expression for a candidate equilibrium is given by
$$
x=P\Psi(x),
$$
for $\psi_k \in \{\psi_k^-,\psi_k^+\}$.
There are $2^n$ such candidates.

Let $\sigma$ be a sign vector $\sigma(k) \in \{-,+\}$ and define the iteration
\begin{equation}\label{itera}
\sigma(k+1) =\mbox{sign}\left [P \Psi(\sigma_k) \right],
\end{equation}
and consider a fixed point of this iteration (if any)
\begin{equation}\label{fixed}
\bar \sigma =\mbox{sign}\left [P \Psi(\bar \sigma) \right].
\end{equation}
The vector $x$ is a rest point if and only if $\sigma=\mbox{sign}(x)$ satisfies \eqref{fixed}. In other words, equation~\eqref{fixed} characterises all the rest points and finding such rest points is equivalent to finding fixed points of the sign iteration. 

The next result follows immediately from the monotone nature of our system, which builds on the condition that $\Psi(y) \ge \Psi(x)$ if $y \ge x$.

\begin{lemma}
Iteration \eqref{itera} is monotone: if $\sigma^A (0) \leq \sigma^B(0)$ are initial sign vectors, then the corresponding iteration satisfies $\sigma^A (k) \leq \sigma^B(k)$.
\end{lemma}


To compute the worst case rest point we initialize $\sigma(0) =[- - \dots -]^\top$. If $\sigma(1)$ has all $-$ signs we have a rest point (all organizations fail). Conversely, let us assume there are $+$ signs. These are nodes that cannot fail. For instance,
$$
\sigma(0)= [- - - - - - - -]^\top
~~\Rightarrow 
\sigma(1)= [- + - - + + - -]^\top
$$
means that nodes $2$, $5$, and $6$ cannot be negative at the equilibrium, due to the monotonicity. These are safe nodes. We denote by $I_{safe}(k)$ the set of safe nodes, namely those that have $+$ signs at the $k$th iteration. Initially, $I_{safe}(0)=\emptyset$, then it cannot reduce to
$$
I_{safe}(0) \subseteq I_{safe}(1) \subseteq I_{safe}(2) \dots
$$ 
Therefore, the iteration converges, stopping when $\sigma(k+1)=\sigma(k)=\sigma^W$. In turn, this means $I_{safe}(k+1) = I_{safe}(k)=I_{safe}^W$. This corresponds to the worst equilibrium $\bar x^W := P \Psi(\sigma^W))$.

By symmetry, we can iterate from the $+$ equilibrium. In this case the safe node set cannot grow, namely,
$$
I_{safe}(0) \supseteq I_{safe}(1) \supseteq I_{safe}(2) \dots
$$ 
The iteration converges to some $\sigma^B$ and the set of safe nodes converges to $I_{safe}^B$. This corresponds to the best equilibrium $\bar x^B := P \Psi(\sigma^B))$. The next result follows immediately from the above.

\begin{lemma}
Any trajectory $\sigma(k)$ satisfies the property
$$
\sigma^W \leq \sigma(k) \leq \sigma^B.
$$
The above bounds hold true also for the fixed points, i.e.,
$$
\sigma^W \leq \bar \sigma \leq \sigma^B.
$$
\end{lemma}

This means that any fixed point for the system satisfies
$$
P \Psi(\sigma^w)) = \bar x^W \leq \bar x \leq \bar x^B = P \Psi(\sigma^B).
$$


If $(P\psi^-)_i > 0$, the $i$th component is always positive (fixed point of the iterative scheme). Likewise, if $(P\psi^+)_i < 0$, the $i$th component is always negative.

\medskip
Consider the trajectory of system~\eqref{eq:sysxr} starting from the negative candidate equilibrium as:

$$x^W(t),~x^W(0) = P \Psi^-.$$

This sequence is monotonically nondecreasing. Indeed,
\ls{\begin{align}
\nonumber x^W(1) & = C x^W(0) + \Psi^- \\
\nonumber       & = P \Psi^- + \Psi(x^W(0))-(I-C)P \Psi^- \\
\nonumber       & = x^W(0) + \Psi(x^W(0))-\Psi^- \geq x^W(0).
\end{align}}
Then, recursively, by monotonicity, we have 
$$
x^W(t+1) \geq x^W(t).
$$
Therefore $x^W(t)$ converges to an equilibrium  $\bar x^W$ from below. Conversely, consider the trajectory of system~\eqref{eq:sysxr} starting from the positive candidate equilibrium as:

$$x^B(t),~x^B(0) = P \Psi^+.$$

This sequence is monotonically nonincreasing, and symmetrically to the above $x^B(t)$ converges to an equilibrium $\bar x^B$ from above. 

Necessarily, these equilibria are related to the bounds $\sigma^W$ and $\sigma^B$
introduced before, then we have
$$
\bar x^W\geq  P\Psi(\sigma^w),~~~~~~\bar x^B \leq P \Psi(\sigma^B),
$$
since these are conditions that hold true for all equilibria.

In fact, the inequalities are satisfied with equal sign. Indeed the initial conditions satisfy
$$x^W(0) = P \Psi^- \leq \bar x^W, \quad x^B(0) = P \Psi^+ \geq \bar x^B,$$
so $x^W(t)$ cannot become greater than $\bar x^W$ and $x^B(t)$ cannot become smaller than $\bar x^B$.

\textit{Remark}. Equilibria $\bar x^W$ and $\bar x^B$ are attractors w.r.t. the initial conditions in othant $2^{n-1}$ and $0$, respectively.


As it is clear from the previous derivation, we further remark that there are points where the $+$ and $-$ are fixed from initialization. In particular, the indices where
$$
(Pr)_i < 0
$$
holds true are $-$ in all iterations. Vice versa, the indices where
$$
(P(r - \beta))_i > 0
$$
holds true are $+$ in all iterations.

\begin{theorem}\label{th:ith} 
Consider system~\eqref{eq:x}. Let $i = 1, \dots, n$.
\begin{itemize}
\item \textbf{Case 1}. Let $(P r)_i < 0$. Then, $\bar x_i < 0$.
\item \textbf{Case 2}. Let $(P(r - \beta))_i \ge 0$. Then, $\bar x_i \ge 0$.
\end{itemize}
\end{theorem}
\textit{Proof}. The proof addresses the above two points one by one. 
\begin{itemize}
\item \textbf{Case 1}. The following 
$$
\bar x_i = \sum_j P_{ij}r_j - \Bigg( \sum_j P_{ij}\beta_j + \cdots \Bigg)
$$ 
is always negative as the first sum is negative and the quantity after the subtraction is positive.
\item \textbf{Case 2}. The following 
$$
\bar x_i = \Bigg( \sum_j P_{ij}r_j - \sum_j P_{ij}\beta_j \Bigg) + \cdots
$$ 
is always positive as the first components in parentheses are positive and the other quantities are also positive.
\end{itemize}
This concludes our proof. \hfill $\blacksquare$

\textit{Remark}. The convergence of the trajectory to a specific configuration of signs means that there exist no oscillations for the dynamical system in the corresponding orthant and the market values converge to the equilibrium point in that orthant.

A direct consequence of Theorem~\ref{th:ith} is the following result, which provides a bound on the number of failed organizations (and saved ones).

\begin{corollary}\label{cor.1}
The number of failed organisations $n_F$ is such that $\mathbb 1^\top_n \phi((I_n - C)^{-1}r) \le n_F \le \mathbb 1^\top_n \phi((I_n - C)^{-1}(r - \beta))$. \hfill $\square$
\end{corollary}
\textit{Proof}. From Theorem~\ref{th.negpos}, $\bar x^{\rm max} = (I_n - C)^{-1}r$ and $\bar x^{\rm min} = (I_n - C)^{-1}(r-\beta)$ such that a generic equilibrium $\bar x$, it holds $\bar x^{\rm min} \le \bar x \le \bar x^{\rm max}$.
Since $\mathbb 1^\top_n \phi(\bar x) = n_F$, the number of failed organizations obeys the stated inequality, equivalent to 
$$\mathbb 1^\top_n \phi(-\ubar V + (I_n - C)^{-1}Dp) \le n_F \le \mathbb 1^\top_n \phi(-\ubar V + (I_n - C)^{-1}(Dp - \beta)).$$
This concludes our proof. \hfill $\blacksquare$ \\

\begin{example}\label{ex4}
Before concluding the paper, we provide one last example in the spirit of \cite{Birge_2021, Elliott_2014}. We now consider system~(\ref{eq:model}) with $N = 9$ organisations, $M=9$ assets. In particular, our analysis involves the cross-holdings among nine countries, i.e., France (FR), Germany (DE), Greece (GR), Italy (IT), Japan (JP), Portugal (PT), Spain ( ES), United Kingdom (GB) and USA (US). 

The matrix of cross-holdings $C$ is summarised in Table~\ref{t:cross}. We assume that $D = I_N$, and $p$ is proportional to the countries GDP as shown in Table~\ref{t:pvalues}. The initial condition $V(0)$ is set to be $V(0) = [15.2838, 19.9137, 0.9863, 9.0642, 28.3350, 0.7829, 8.8020,$ $12.1361, 59.8130]^\top$. We set $\beta = 0.5 \; \mathbb 1_{20}$ and $\ubar V = 10 \; \mathbb 1_{20}$.
\begin{table} [h]
\caption{Table providing the values of the matrix of cross-holdings $C$, adapted from \cite{Birge_2021}.}
\begin{center}
   \begin{tabular}{|c|c|c|c|c|c|c|c|c|c|}
  \hline
 & FR & DE & GR & IT & JP & PT & ES & GB & US \\ \hline \hline
   FR & $0$ & $.03$ & $.01$ & $.07$ & $.01$ & $.04$ & $.04$ & $.05$ & $.04$ \\ \hline
   DE & $.04$ & $0$ & $.06$ & $.03$ & $.00$ & $.05$ & $.04$ & $.09$ & $.04$ \\ \hline
   GR & $.00$ & $.00$ & $0$ & $.00$ & $.00$ & $.00$ & $.00$ & $.00$ & $.00$ \\ \hline
   IT & $.01$ & $.03$ & $.00$ & $0$ & $.00$ & $.01$ & $.02$ & $.01$ & $.00$ \\ \hline
   JP & $.04$ & $.02$ & $.00$ & $.02$ & $0$ & $.01$ & $.01$ & $.06$ & $.10$ \\ \hline
   PT & $.00$ & $.00$ & $.00$ & $.00$ & $.00$ & $0$ & $.00$ & $.00$ & $.00$ \\ \hline
   ES & $.01$ & $.02$ & $.01$ & $.02$ & $.00$ & $.15$ & $0$ & $.09$ & $.02$ \\ \hline
   GB & $.03$ & $.02$ & $.01$ & $.01$ & $.01$ & $.02$ & $.01$ & $0$ & $.04$ \\ \hline
   US & $.04$ & $.02$ & $.01$ & $.02$ & $.02$ & $.02$ & $.02$ & $.09$ & $0$ \\ \hline
    \end{tabular}\\
  \end{center}
  \label{t:cross}
\end{table}

\begin{table} [h]
\caption{Original primitive asset values $p$ \cite{Birge_2021}.}
\begin{center}
   \begin{tabular}{|c|c|}
  \hline
 Country & Relative GDP \\ \hline \hline
   FR & $12.29$  \\ \hline
   DE & $16.81$  \\ \hline
   GR & $1.02$  \\ \hline
   IT & $9.30$  \\ \hline
   JP & $20.00$  \\ \hline
   PT & $1.00$  \\ \hline
   ES & $6.00$  \\ \hline
   GB & $12.99$  \\ \hline
   US & $75.70$  \\ \hline
    \end{tabular}\\
  \end{center}
  \label{t:pvalues}
\end{table}

We show the behaviour of the nine countries and their convergence to $\overline V \ge 0$. This is in accordance with Theorem~\ref{th.1}. Figure~\ref{fig:countries} shows this scenario.
\begin{figure}[t]
    \centering
    \includegraphics[width=0.45\textwidth]{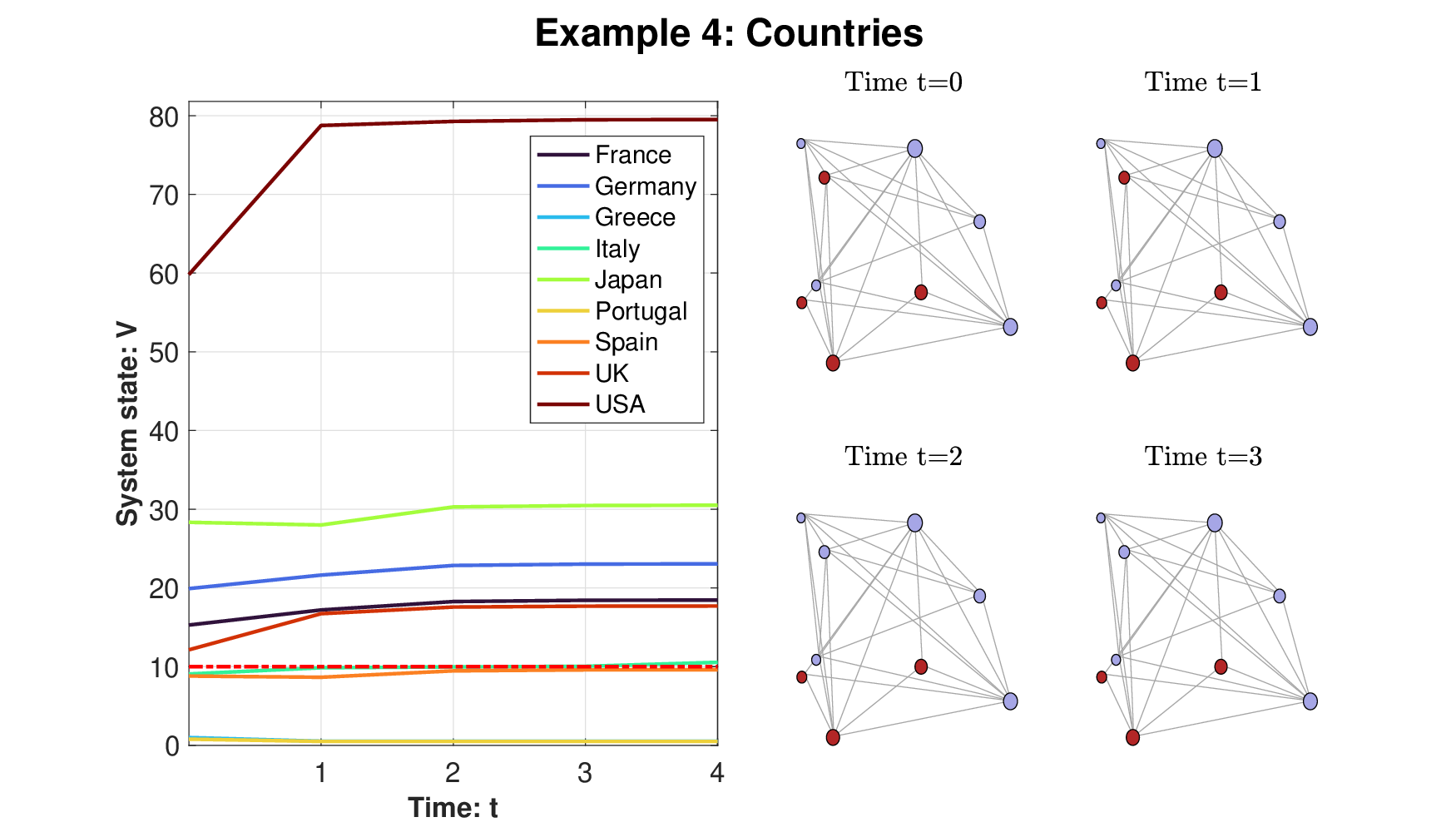}
    \caption{Behaviour of the system representing 9 countries and the cross-holdings among them.}
    \label{fig:countries}
\end{figure}

\ls{\textit{Remark}. This last example provides an application of Theorem~\ref{th.1} in a real life scenario where the organizations are represented by a set of countries. Analogously to Example~\ref{ex2}, we can determine the evolution of the equity values since the condition of Theorem~\ref{th.1} holds true, meaning that the equity value of every country will remain positive at all times.}
\end{example}

\section{Conclusions}\label{sec:conc}
In this paper, we study the propagation of failures in financial systems. \ls{Future works include: i) the characterization of the invariance of each orthant of the $2^n$ space and of the equilibria in each orthant, ii) the study of the worst-case scenario where all organizations fail and the conditions to prevent it, and iii) asset investments as feedback control design.}


\section*{Acknowledgments}
LS has been partly supported by the HUMAT Research Project, financed by the US Army Research Lab, USA, in collaboration with the Alan Turing Institute, UK.
DB has been supported the SMiLES Research Project, part of the Research Programme Sustainable Living Labs, which is co-financed by the Dutch Research Council (NWO), the Ministry of Infrastructure and Water Management, The Netherlands, the Taskforce for Applied Research (SIA), The Netherlands, and the Top Sector Logistics, The Netherlands.
FB and PC has been partly supported by the Italian grant PRIN 2017 ``Monitoring and Control Underpinning the Energy-Aware Factory of the Future: Novel Methodologies and Industrial Validation'' (ID 2017YKXYXJ).
FB has been supported by the European Union - NextGenerationEU.




\begin{thebibliography}{00}



\bibitem{Ahn_2019} D. Ahn, K. Kim, Optimal intervention under stress scenarios: A case of the Korean financial system, Oper. Res. Lett. 47 (4) (2019) 257--263.

\bibitem{Acemoglu_2012} D. Acemoglu, V.M. Carvalho, A. Ozdaglar and A. Tahbaz-Salehi, The network origins of aggregate fluctuations, Econometrica 80 (2012) 1977--2016.

\bibitem{Acemoglu_2015} D. Acemoglu, A. Ozdaglar and A. Tahbaz-Salehi, Systemic risk and stability in financial networks, Am. Econ. Rev. 105 (2) (2015) 564--608.

\bibitem{Birge_2021} J.R. Birge, Modeling Investment Behavior and Risk Propagation in Financial Networks, SSRN (2021).

\bibitem{Brioschi_1989} F. Brioschi, L. Buzzacchi and M.G. Colombo, Risk capital financing and the separation of ownership and control in business groups, J. Bank. Financ. 13 (4) (1989) 747--772.

\bibitem{Cabrales_2017} A. Cabrales, P. Gottardi and F. Vega-Redondo, Risk sharing and contagion in networks, Rev. Financ. Stud. 30 (9) (2017) 3086--3127.

\bibitem{Calafiore_2022} G.C. Calafiore, G. Fracastoro and A.V. Proskurnikov, Control of dynamic financial networks, IEEE Control Syst. Lett. 6 (2022) 3206--3211.

\bibitem{Carvalho_2008} V.M. Carvalho, Aggregate fluctuations and the network structure of intersectoral trade, The University of Chicago.

\bibitem{Chen_2021} H. Chen, T. Wang and D.D. Yao, Financial network and systemic risk -- a dynamic model, Prod. Oper. Manag. 30 (8) (2021) 2441--2466.

\bibitem{Eisenberg_2001} L. Eisenberg and T.H. Noe, Systemic risk in financial systems, J. Manag. Sci. 47 (2) (2001) 236--249.

\bibitem{Elliott_2014} M. Elliott, B. Golub and M.O. Jackson, Financial networks and contagion, Am. Econ. Rev. 104 (10) (2014) 3115--3153.

\bibitem{Elsinger_2009} H. Elsinger, Financial networks, cross holdings, and limited liability, Working Papers 156, \"{O}sterreichische Nationalbank (Austrian Central Bank) (2009).

\bibitem{Farina_2000} L. Farina and S. Rinaldi, Positive linear systems -- Theory and applications, Wiley \& Sons, New York, 2000.

\bibitem{Fischer_2014} T. Fischer, No-arbitrage pricing under systemic risk: Accounting for cross-ownership, Math. Financ. 24 (1) (2014) 97--124.

\bibitem{Glasserman_2015} P. Glasserman and H.P. Young, How likely is contagion in financial networks? J. Bank. Financ. 50 (2015) 383--399.

\bibitem{Glasserman_2016} P. Glasserman and H.P. Young, Contagion in financial networks, J. Econ. Lit. 54 (3) (2016) 779--831.

\bibitem{Karl_2014} S. Karl and T. Fischer, Cross-ownership as a structural explanation for over- and underestimation of default probability. Quant. Finance 14 (6) (2014) 1031--1046.

\bibitem{Ramirez_2022} S. Ramirez, M.v.d. Hoven and D. Bauso, A stochastic model for cascading failures in financial networks, IEEE Trans. Control. Netw. Syst. doi: 10.1109/TCNS.2023.3256273.

\bibitem{Rogers_2013} L.C.G. Rogers and L.A.M. Veraart, Failure and rescue in an interbank network, J. Manag. Sci. 59 (4) (2013) 882--898.

\bibitem{Weber_2017} S. Weber and K. Weske, The joint impact of bankruptcy costs, fire sales and cross-holdings on systemic risk in financial networks, Probability, Uncertainty and Quantitative Risk 2 (9) (2017).


%






\end{thebibliography}


\end{document}